\documentclass[a4 paper, 12pt]{amsart}
\usepackage{amsmath,amsfonts,amssymb,amsthm,multicol,enumerate,hyperref,bookmark,graphicx,float,tikz, footmisc, fancyhdr}
\usepackage[top=2in, bottom=1.5in, left=1in, right=1in]{geometry}
\newtheorem{thm}{Theorem}[section]

\numberwithin{equation}{section}
\hypersetup{colorlinks=true, citecolor= green, linkcolor=blue,urlcolor=red}

\begin{document}
\title[Finite Time Stability Analysis]
{Finite Time Stability Analysis of Non-Linear Fractional Order With Multi State Time Delay}
\author[Supriyo Dutta and N. Sukavanam]{}
\maketitle
\begin{center}
{\bf Supriyo Dutta and N. Sukavanam} \\
Department of mathematics\\Indian Institute of Technology Roorkee\\ dutta.1@iitj.ac.in, nsukvfma@iitr.ernet.in
\end{center}

\begin{abstract}
Sufficient condition for the stability of a fractional order semi-linear system with multi-time delay is proposed.
\end{abstract}

\section{Introduction}

Stability of linear and non-linear control system is an important area of research in control theory. Stability of finite dimensional system of ordinary and fractional differential equation is discussed in \cite{Robert2004} and \cite{Ivo2009}. Gronwall Inequality was generalized in \cite{Ye2007} and is used to discuss the finite time stability in \cite{Liu2011}. In \cite{Harry1955} Mittage-Leffler function and its properties are discussed in details.

The main result of this paper provides a condition for finite time stability of a fractional order system with multi-state time delay. It depends on the result of Generalized Gronwall' inequality for fractional order differential equation.

\section{Fundamental concepts}

In this section we provide an overviews on the fundamentals related to Riemann-Liouvill fractional integral and derivatives \cite{RossBertram1993}, fractional order differential equations and fractional order system, generalization of the Granwall' inequality for fractional order differential equation \cite{Ye2007} and finite time stability.

\subsection{Fractional Integral and Differential}
Riemman Liouville Fractional integral of a function $f :\mathbb{R}\rightarrow\mathbb{R}$ of order $\nu$ is defined as \cite{RossBertram1993} \cite{Antoly2006}
$$_0\mathrm{D}_{t}^{-\nu}f(t)=\frac{1}{\Gamma(\nu)}\int_{0}^{t}(t-\xi)^{\nu-1}f(\xi)\mathrm{d}\xi\;;\;Re(\nu)>0$$
The $\mu$-th order Riemman Liouville fractional derivative is given by
$$\mathrm{D}^\mu f(t)=\left\{\begin{array}{c c}\frac{1}{\Gamma(\nu)}\frac{\mathrm{d}^m}{\mathrm{d}t^m}\int_0^t(t-\xi)^{\nu-1}f(\xi) \mathrm{d}\xi&\;\mbox{when $m-1<\mu < m$}\\\frac{\mathrm{d}^m}{\mathrm{d}t^m}f(t)\;&\mbox{when $\mu = m \in \mathbb{N} $ } \end{array} \right. $$
Consider the semi-linear fractional order differential equations
\begin{equation}\label{m1}
\begin{split}
\mathrm{D}^q x(t) &= A_0x(t)+\sum_{i=1}^pA_ix(t-\tau_i)+B_0u(t)+f(t,x) ; \; t\geq 0 \\
x(t) &= \psi_x(t) ; \; t \in [-\tau,0]
\end{split}
\end{equation}
on the space $(\mathrm{C}[-\tau,T],\mathbb{R}^n)$ with the uniform norm defined as   $$\|x\|_=\max\{|x_1(t)|,|x_2(t)|,\dots,|x_n(t)|\}.$$
Where
\begin{gather*}
\begin{split}
& \tau = max\{\tau_1,\tau_2,\dots,\tau_p\} ; \;\text{$\tau_i > 0$ are contants.}\\
& x(t) :\mathbb{R} \rightarrow \mathbb{R}^n; \; \text{is a $n\times 1 $ vector and}\; x \in C[[-\tau, T] : \mathbb{R}^n]\\
& \mathrm{D}^qx(t) = (\mathrm{D}^{q}x_1(t), \mathrm{D}^{q}x_2(t), \dots, \mathrm{D}^{q}x_n(t))^t ; \; \text{is an $n \times 1$ vector.}\\
& (A_i)_{n \times n} = \{a_{j , k}\}_i;\; \text{are constant matrix for}\;i=0, 1, \dots, p.\\
\end{split}
\end{gather*}
$f(t,x)_{n \times 1} : \mathbb{R} \times \mathbb{R}^n \rightarrow \mathbb{R}^n$; satisfies Lipchitz condition on w.r.t $x$.
$$\therefore \|f(t,x) - f(t,y)\| \le L \|x - y \|$$
$$\Rightarrow \|f(t,x)\| \le L\|x\| + m , $$
where $ m = \|f(t, \theta)\|, \theta$ is null vector.

\subsection{Generalized Gronwall Inequality}\cite{Ye2007}
Suppose $y(t):\mathbb{R}\rightarrow\mathbb{R}$ and $a(t):\mathbb{R}\rightarrow\mathbb{R}$ are non-negative and integrable in every closed and bounded subinterval of $[0,T)$ and $g(t):[0,T)\rightarrow\mathbb{R}$ is non-negative, non-decreasing, continuous and bounded function such that
$$y(t)\leq a(t) + g(t)\int_0^t(t-s)^{q-1}y(s)\mathrm{d}s$$
then for $0 \le t \le T$
$$y(t)\leq a(t)+\int_0^t\left[\sum_{n=1}^\infty\frac{g(t)\Gamma(q)}{\Gamma(nq)}(t-s)^{nq-1}a(s)\right]\mathrm{d}s$$
Moreover if $a(t)$ is nondecreasing then
$$y(t) \leq a(t)E_q(g(t).\Gamma(q)t^q).$$
Where $E_q$ is the Mittag-Leffler function, defined by
$$E_q(z)=\sum_{k=0}^\infty\frac{z^k}{\Gamma(kq+1)}.$$

\subsection{Finite Time Stability}\cite{Liu2011}
The system given by \eqref{m1} is finite time stable w.r.t. $\{\delta, \varepsilon, q_u, J\}$, if $\|\psi_x\|<\delta$ and $\|u(t)\|<q_u, \forall t \in J = [0,T]$ imply $\|x(t)\|<\varepsilon, \forall t\in J$. \\

The system \eqref{m1} is called homogeneous if $u(t)=0$. In this case the system will be finite time stable if $\|\psi_x\|<\delta$ imply $\|x(t)\|<\varepsilon, \forall t \in J$.\\

In \cite{Liu2011} it has been shown that the linear system \eqref{m1}is stable if the following condition is satisfied. \footnote{The condition derived in this paper is stronger than that of \cite{Liu2011}}
$$\left[1 + \frac{(n+1)\sigma t^q}{\Gamma(q+1)} + \frac{q_u b_0 t^q}{\delta \Gamma(q+1)}\right]E_q((n+1)\sigma t^q) < \frac{\varepsilon}{\delta}$$

\section{Main Result}
\begin{thm}\label{m2}
The system given by \eqref{m1} is finite time stable w.r.t. $\{ \delta , \varepsilon, q_u , J \}$, if the following condition is satisfied
$$\left[1+ \frac{( m + b q_u )T^q}{\delta \Gamma( q + 1 )}\right] E_q\{(L+\sigma(p+1))T^q\} < \frac{\varepsilon}{\delta}$$
\end{thm}
\begin{proof}
The solution of \eqref{m1} is given by
\begin{equation*}
\begin{split}
x(t) =& x(0) + \frac{1}{\Gamma(q)}\int_0^t (t-s)^{q-1} \left\{A_0 x(s) + \sum_{i=1}^p A_i x(s-\tau_i)+ B_0u(s)+ f(s,x)\right\}\mathrm{d}s\\
\Rightarrow \|x(t)\| \leq& \|x(0)\| + \frac{1}{\Gamma(q)}\int_0^t (t-s)^{q-1} \{ \|A_0\|\| x(s)\| + \sum_{i=1}^p \|A_i\|\| x(s-\tau_i)\|\\   &+ \|B_0\|\|u(s)\| + \|f(s,x)\| \} \mathrm{d}s
\end{split}
\end{equation*}
Now let $\sigma_{max}(A)$ is the largest singular value of the matrix $A$ and
\begin{gather*}
\begin{split}
\sigma &= \max _{0\leq i \leq p}\{\sigma_{max}(A_i)\} \\
b_0 &= \sigma_{max}(B_0)\\
\therefore \|A_i\| &\leq \sigma ; \; \forall i = 1,2,\dots,p. \\
\end{split}
\end{gather*}
\begin{equation*}
\begin{split}
\| x( t ) \| \leq& \| x( 0 ) \| + \frac{1}{\Gamma( q )} \int_0^t ( t - s )^{q - 1} \{ \sigma \| x( s ) \| + \sum_{i=1}^p \sigma \| x ( s - \tau_i ) \|\\
& + b_0 \| u( s ) \| + \| f( s , x ) \| \} \mathrm{d}s \\
\leq & \| \psi(0) \| + \frac{1}{\Gamma(q)} \int_0^t ( t - s )^{q - 1} \{ \sigma \| x(s) \| + \sum_{i=1}^p \sigma  \| x(s) \| \\
& + b_0 \| u( s ) \| + m + L \| x( s ) \| \} \mathrm{d}s \\
\leq & \| \psi( 0 ) \| + \frac{1}{\Gamma(q)} \int_0^t (t - s)^{q - 1}\{\sigma ( p + 1 )\|x( s ) \|\\
& + b_0 q_u + m + L \| x( s ) \| \} \mathrm{d}s \\
\leq & \| \psi( 0 ) \| + \frac{ m + b_0 q_u }{\Gamma( q )} \int_0^t ( t - s )^{q - 1} \mathrm{d}s +\\
& \frac{ L + \sigma ( p + 1 )}{\Gamma( q )} \int_0^t ( t - s )^{q - 1} \| x( s ) \| \mathrm{d}s \\
\leq & \left\{ \delta + \frac{( m + b_0 q_u ) t^q}{\Gamma( q + 1 )}\right\} + \frac{ L + \sigma (p+1)}{\Gamma( q )} \int_0^t ( t - s )^{q - 1} \| x( s ) \| \mathrm{d}s \\
\| x( t ) \| \leq& a(t) + g(t) \int_0^t ( t - s )^{ q - 1 } \| x( s ) \| \mathrm{d}s\\
\text{Where} \; a( t ) =& \delta + \frac{( m + b_0 q_u )t^q }{\Gamma( q + 1 )} \;  \\
\text{and} \; g( t ) =& \frac{ L + \sigma ( p + 1 )}{\Gamma( q )}
\end{split}
\end{equation*}
Now $a(t)$ is nonnegative and locally integrable. $g(t)$ is non negative bounded and nondecreasing. Now by the Gronwall' inequality we shall say that
\begin{equation*}
\begin{split}
\|x(t)\| \leq & a(t)E_q(g(t)\Gamma(q)t^q)\\
= & \left\{\delta + \frac{(m + b_0q_u ) t^q } { \Gamma ( q + 1 ) } \right \} E_q \left\{ \frac{L + \sigma ( p + 1 ) }{ \Gamma ( q ) } \Gamma(q)t^q\right\}\\
\leq& \delta\left\{1+\frac{(m+ b_0 q_u)T^q}{\delta\Gamma(q+1)}\right\}E_q\left\{(L+\sigma(p+1))T^q\right\}\\
\end{split}
\end{equation*}
So for the finite time stability we require $\|x(t)\|<\varepsilon$ if
$$\left\{1+\frac{(m+b_0q_u)T^q}{\delta\Gamma(q+1)}\right\}E_q\{(L+\sigma(p+1))T^q\} < \frac{\varepsilon}{\delta}$$
\end{proof}
\textbf{Special Cases}
\begin{enumerate}
\item
If the control term $u(t)=0$, then $q_u=0$. Then theorem \eqref{m2} will take the following form
$$\left\{1+\frac{m}{\delta \Gamma(q+1)}T^q\right\}E_q\left\{(L+\sigma(p+1))T^q\right\} \leq \frac{\varepsilon}{\delta}$$
\item
If $f(t,\theta)=0$ i.e. $ m = 0 $ them theorem \eqref{m2} will take the following form
$$\left\{1+\frac{b_0q_u}{\delta \Gamma(q+1)}T^q\right\}E_q\left\{(L+\sigma(p+1))\right\} \leq \frac{\varepsilon}{\delta}$$
\item
Let there is no non-linear term i.e. $f(t)=0$. So $L=0$ and $ m = 0 $. Then the theorem \eqref{m2} will take the following form
$$\left\{1+\frac{b_0q_u}{\delta \Gamma(q+1)}T^q\right\}E_q\left\{\sigma(p+1)T^q\right\} \leq \frac{\varepsilon}{\delta}$$
\end{enumerate}

\bibliographystyle{plain}

\end{document}